\documentclass[11pt]{article}

\usepackage[letterpaper,margin=1.1in]{geometry}

\usepackage{natbib}
\usepackage{comment}
\usepackage{amsfonts}
\usepackage{amsthm}
\usepackage{amssymb}
                 
\newcommand{\R}{\mathbb{R}}

\newcommand{\Exp}{\mathbb{E}}

\newcommand{\vx}{\bm{x}} 
\newcommand{\vX}{\bm{X}}
\newcommand{\vy}{\bm{y}} 

\newcommand{\vp}{\bm{p}}
                 
\usepackage{bm}

\newcommand{\vxi}{\bm{\xi}}

\newcommand*{\SPSP}{\ensuremath{\left(\mathbf{SPSP}\right)}}

\newcommand*{\SPSPM}{\ensuremath{\left(\mathbf{SPSP-M}\right)}}

\newcommand*{\MP}{\ensuremath{\left(\mathbf{MP}\right)}}
\newcommand{\cvar}{\operatorname{CVaR}}
                 
\newcommand{\ds}{\displaystyle}
\newcommand{\vsp}{\vspace{0.5ex}}

\usepackage{mathtools}
\usepackage{multicol}
\usepackage{longtable}
\usepackage{booktabs} 

\newtheorem{theorem}{Theorem}
\newtheorem{lemma}{Lemma}
\newtheorem{corollary}{Corollary}

\newcommand{\tilim}{\color{blue} $^\dagger$}

\title{\bf Stochastic Planning and Scheduling \\ with Logic-Based Benders Decomposition}

\author{\"{O}zg\"{u}n El\c{c}i and J. N. Hooker\\
\ \\
{\normalsize Carnegie Mellon University} \\
{\normalsize elci@cmu.edu, jh38@andrew.cmu.edu}
} 

\date{December 2020}

\begin{document}
\maketitle

\begin{abstract}
	We apply logic-based Benders decomposition (LBBD) to two-stage stochastic planning and scheduling problems in which the second-stage is a scheduling task.   We solve the master problem with mixed integer/linear programming and the subproblem with constraint programming.  As Benders cuts, we use simple nogood cuts as well as analytical logic-based cuts we develop for this application.   We find that LBBD is computationally superior to the integer L-shaped method, with a branch-and-check variant of LBBD faster by several orders of magnitude, allowing significantly larger instances to be solved.  This is due primarily to computational overhead incurred by the integer L-shaped method while generating classical Benders cuts from a continuous relaxation of an integer programming subproblem.  To our knowledge, this is the first application of LBBD to two-stage stochastic optimization with a scheduling second-stage problem, and the first comparison of LBBD with the integer \mbox{L-shaped} method.  The results suggest that LBBD could be a promising approach to other stochastic and robust optimization problems with integer or combinatorial recourse. 
\end{abstract}




%


\section{Introduction}
\label{Sec:Intro}

Benders decomposition has seen many successful applications to two-stage stochastic optimization, where it typically takes the form of the {\em L-shaped method} \citep{benders1962partitioning,van1969shaped}.  It offers the advantage that the second-stage problem decouples into a separate problem for each possible scenario, allowing much faster computation of the recourse decision.  

A limitation of classical Benders decomposition, however, is that the subproblem must be a linear programming problem, or a continuous  nonlinear programming problem in the case of ``generalized'' Benders  decomposition \citep{geoffrion1972generalized}.  This is necessary because the Benders cuts are derived from dual multipliers (or Lagrange multipliers) in the subproblem.  Yet in many problems, the recourse decision is a combinatorial optimization problem that does not yield dual multipliers.  This issue has been addressed by the {\em integer \mbox{L-shaped} method} \citep{laporte1993integer}, which formulates the subproblem as a mixed integer/linear programming (MILP) problem and obtains dual multipliers from its linear programming (LP) relaxation.  To ensure finite convergence, classical Benders cuts from the LP relaxation are augmented with ``integer cuts'' that simply exclude the master problem solutions enumerated so far.  

Unfortunately, a combinatorial subproblem may be difficult to model as an MILP, in the sense that many variables are required, or the LP relaxation is weak.  This is particularly the case when the recourse decision poses a scheduling problem.  We therefore investigate the option of applying {\em logic-based Benders decomposition} (LBBD) to problems with a second-stage scheduling decision \citep{Hoo00,hooker2003logic}, because it does not require dual multipliers to obtain Benders cuts.  Rather, the cuts are obtained from an ``inference dual'' that is based on a structural analysis of the subproblem.  This allows the subproblem to be solved by a method that is best suited to compute optimal schedules, without having to reformulate it as an MILP.  

We investigate the LBBD option by observing its behavior on a generic planning and scheduling problem in which scheduling takes place after the random events have been observed.  The planning element is an assignment of tasks to facilities that occurs in the first-stage. Tasks assigned to each facility are then scheduled in the second-stage subject to time windows.  We assume that the task processing time is a random variable, but the LBBD approach is easily modified to accommodate other random elements, such as the release time.  The subproblem decouples into a separate scheduling problem for each facility and each scenario.  For greater generality, we suppose the recourse decision is a {\em cumulative} scheduling problem in which multiple tasks can run simultaneously on a single facility, subject to a limit on total resource consumption at any one time.  

We solve the first-stage problem by MILP, which is well suited for assignment problems.  More relevant to the present study is our choice to solve the scheduling subproblem by constraint programming (CP), which has proved to be effective on a variety of scheduling problems, perhaps the state of the art in many cases.  We therefore formulate the subproblem in a CP modeling language rather than as an MILP.  In view of the past success of LBBD on a number of deterministic planning and scheduling problems, we test the hypothesis that it can obtain similar success on stochastic problems with many scenarios.  Our computational study focuses on the minimum makespan problem as a proof of concept, but we show how LBBD is readily modified to accommodate other objectives, such as minimizing total tardiness or total assignment cost.  We also derive new logic-based Benders cuts for the minimum makespan problem that have not been used in previous work.  

In addition to standard LBBD, we experiment with {\em branch and check}, a variation of LBBD that solves the master problem only once and generates Benders cuts on the fly during the MILP branching process \citep{Hoo00,Tho01}.  We find that both versions of LBBD are superior to the integer L-shaped method.  In particular, branch and check is faster by several orders of magnitude, allowing significantly larger instances to be solved.  We also conduct a variety of tests to identify factors that explain the superior performance of LBBD, the relative effectiveness of various Benders cuts, and the impact of modifying the integer L-shaped method in various ways.  To our knowledge, this is the first computational comparison between LBBD and the integer L-shaped method on any kind of stochastic optimization problem.  It also appears to be the first application of LBBD to two-stage stochastic optimization with a scheduling second-stage problem.

The remainder of this paper is organized as follows.
We introduce the stochastic planning and scheduling problem in
Section \ref{Sec:Problem_definition}. This is followed
by Section \ref{Sec:Algorithm} where we propose the logic-based
Benders decomposition based solution methods for solving
three variants of the stochastic planning and scheduling problem.
We present the computational results in Section
\ref{Sec:CompStudy} and give our concluding remarks in
Section \ref{Sec:Conclusion}.

\section{Previous Work}

A wide range of problems can be formulated as two-stage
stochastic programs. For theory and various applications, 
we refer the reader to \cite{birge2011introduction}, 
\cite{shapiro2009lectures}, \cite{prekopa2013stochastic},
and the references therein.
Allowing discrete decisions in the second-stage problem
significantly expands the applicability of the two-stage
stochastic framework, as for example to last-mile relief network design
\citep{noyan2015stochastic} and
vehicle routing with stochastic travel times
\citep{laporte1992vehicle}.

Benders decomposition \citep{benders1962partitioning} has long been applied to large-scale optimization problems \citep{geoffrion1974multicommodity,cordeau2001benders,binato2001new,contreras2011benders}.   \cite{rahmaniani2017benders} provide an excellent survey of enhancements to the classical method.  In particular, it has been applied to two-stage stochastic programs with linear recourse by means of the L-shaped method \citep{van1969shaped}.  Its applicability was extended to integer recourse by the integer \mbox{L-shaped} method of \cite{laporte1993integer}, which was recently revisited and improved by \cite{angulo2016improving} and
\cite{li2018improved}. Other Benders-type algorithms that have been proposed for integer
recourse include disjunctive decomposition \citep{sen2005c3} and decomposition with
parametric Gomory cuts \citep{gade2014decomposition}. The essence of these two methods is to convexify the integer second-stage problem using disjunctive cuts and Gomory cuts, respectively. Still other decomposition-based methods in the literature include 
progressive hedging for multi-stage stochastic convex programs
\citep{rockafellar1991scenarios} and a dual decomposition method for multi-stage stochastic programs with mixed-integer variables \citep{caroe1999dual}. We refer the reader to \cite{kuccukyavuz2017introduction} for a review of two-stage stochastic mixed-integer programming.

Logic-based Benders decomposition was introduced by \cite{Hoo00} and further developed in \cite{hooker2003logic}.  Branch and check, a variant of LBBD, was also introduced by \cite{Hoo00} and first tested computationally by \cite{Tho01}, who coined the term ``branch and check.''  A general exposition of both standard LBBD and branch and check, with an extensive survey of applications, can be found in \cite{Hoo19a}.  A number of these applications have basically the same mathematical structure as the planning and scheduling problem studied here, albeit generally without a stochastic element.

In more recent work, \cite{atakan2017minimizing} focus on a one-stage stochastic model for single-machine scheduling in which they minimize the value-at-risk of several random performance measures. \cite{Bulbul2014} consider a two-stage chance-constrained mean-risk stochastic programming model for single-machine scheduling problem, but the scheduling decisions do not occur in the second-stage.  Rather, the second-stage problem is a simple optimal timing problem that can be solved very rapidly.  The deterministic version of the planning and scheduling problem we consider here is solved by LBBD in \cite{hooker2007planning} and \cite{CirCobHoo16}.  We rely on some techniques from these studies.  

We are aware of three prior applications of LBBD to stochastic optimization.  \cite{LomMilRugBen10} use LBBD to assign computational tasks to chips and to schedule the tasks assigned to each chip.  However, this is not an instance of two-stage stochastic optimization, because the stochastic element appears in the first-stage assignment problem and is replaced with its deterministic equivalent. \cite{fazel2013solving} solve a stochastic location-routing problem with LBBD, but there is no actual recourse decision in the second-stage, which only penalizes vehicles if the route determined by first-stage decisions exceeds their threshold capacity.  \cite{GuoBodAleUrb19} use LBBD to schedule patients in operating rooms, where the random element is the surgery duration.  Here the scheduling takes place in the master problem, where patients are assigned operating rooms and surgery dates.  The subproblem checks whether there is time during the day to perform all the surgeries assigned to a given operating room, and if not, finds a cost-minimizing selection of surgeries to cancel on that day.  Unstrengthened nogood cuts are used as  LBBD cuts, along with classical Benders cuts obtained from a network flow model of the subproblem that is obtained from a binary decision diagram.

The present study therefore appears to be the first application of LBBD to two-stage stochastic optimization with scheduling in the second-stage.  It is also the first to compare any application of stochastic LBBD with the integer L-shaped method.

\section{The Problem}
\label{Sec:Problem_definition}

We study a two-stage stochastic programming problem that, in general, has the following form:
\begin{equation}
\min_{\vx\in X} \big\{f(\vx) + \Exp_{\omega}[Q(\vx,\omega)] \big\}
\label{TwoStage:General}
\end{equation}
where $Q(\vx,\omega)$ is the optimal value of the second-stage problem:
\begin{equation}
Q(\vx,\omega) = \min_{\vy\in Y(\omega)} \big\{g(\vy)\big\}
\end{equation}
Variable $\vx$ represents the first-stage decisions, while 
$\vy$ represents second-stage decisions that are made after the random variable $\omega$ is realized.  We suppose that $\omega$ ranges over a finite set $\Omega$ of possible scenarios, where each scenario $\omega$ has probability $\pi_{\omega}$.  The first-stage problem (\ref{TwoStage:General}) may therefore be written
\[
\min_{\vx\in X} \Big\{f(\vx) + \sum_{\omega\in \Omega} \pi_{\omega} Q(\vx,\omega) \Big\}
\]

We consider a generic 
planning and scheduling problem in which the first-stage assigns tasks to facilities, and the second-stage schedules the tasks assigned to each facility.  The objective is to minimize makespan or total tardiness.  We assume that only the processing times are random in the second-stage, but a slight modification of the model allows for random release times and/or deadlines as well.  

We therefore suppose that each job $j$ has a processing time $p^{\omega}_{ij}$ on facility $i$ in scenario $\omega$ and must be processed during the interval $[r_j,d_j]$.  For greater generality, we allow for cumulative scheduling, where each job $j$ consumes resources $c_{ij}$ on facility $i$, and the total resource consumption must not exceed $K_i$.

To formulate the problem we let variable $x_j$ be the facility to which job $j\in J$ is assigned.  The first-stage problem is
\begin{equation}
\min_{\vx} \Big\{g(\vx) + \sum_{\omega\in \Omega} \pi_{\omega}Q(\vx,\omega) \; \Big| \; x_j\in I, \;\mbox{all} \; j\in J \Big\}
\label{eq:P&S}
\end{equation}
where $I$ indexes the facilities.  In the second-stage problem, we let $s_j$ be the time at which job $j$ starts processing. 
We also let $J_i(\vx)$ be the set of jobs assigned to facility $i$, so that $J_i(\vx)=\{j\in J \;|\; x_j=i\}$. 
Thus
\[
Q(\vx,\omega) = \min_{\bm{s}} \Big\{ h(\bm{s},\vx,\omega) \;\Big|\; s_j\in [r_j,d_j-p^{\omega}_{x_jj}], \;\mbox{all}\;j\in J; \;\; 
\hspace{-3ex} \sum_{\substack{j\in J_i(\vx)\\0\leq t\leq s_j+p^{\omega}_{ij}}} \hspace{-3ex} c_{ij} \leq K_i, \; \mbox{all}\; i\in I,\;\mbox{all}\; t\Big\}
\]

The two-stage problem (\ref{TwoStage:General}) is risk-neutral in the sense that it is concerned with minimizing expectation.  However, the LBBD approach presented here can be adapted to a more general class of problems that incorporate a dispersion statistic $\mathbb{D}_{\omega}$ that measures risk, such as variance, as in the classical Markovitz model \cite{}.  Then the problem (\ref{TwoStage:General}) becomes
\begin{equation}
\min_{\vx\in X} \big\{f(\vx) + (1-\lambda)\Exp_{\omega}[Q(\vx,\omega)] + \lambda \mathbb{D}_{\omega}[Q(\vx,\omega)]\big\}
\label{TwoStage:Risk}
\end{equation}
and the first-stage planning and scheduling problem (\ref{eq:P&S}) becomes 
\begin{equation}
\min_{\vx} \Big\{g(\vx) + (1-\lambda)\sum_{\omega\in \Omega} \pi_{\omega}Q(\vx,\omega) + \lambda\mathbb{D}_{\omega}\big[Q(\vx,\omega)\big] \; \Big| \; x_j\in I, \;\mbox{all} \; j\in J \Big\}
\label{eq:P&Srisk}
\end{equation}
Formulations (\ref{TwoStage:Risk}) and (\ref{eq:P&Srisk}) also accommodate robust optimization, as for example when $\lambda=1$ and 
\[
\mathbb{D}_{\omega}(Q(\vx,\omega)) = \max_{\omega\in\Omega}\{Q(\vx,\omega)\}
\] 
and $\Omega$ is an uncertainty set.   See \cite{ahmed2006convexity} for a discussion of various tractable and intractable risk measures.

\section{Logic-based Benders Decomposition}
\label{Sec:Algorithm}

Logic-based Benders decomposition (LBBD) is designed for problems of the form
\begin{equation}
\label{LBBD:general}
\min_{\vx,\vy} \big\{ f(\vx,\vy) \; \big| \; C(\vx,\vy), \; \vx\in D_x, \; \vy\in D_y \big\}
\end{equation}
where $C(\vx,\vy)$ denotes a set of constraints
that contain variables $\vx$ and $\vy$, and 
$D_y$ and $D_x$ represent variable domains. The rationale behind
dividing the variables into two groups is that
once some of the decisions are fixed by setting $\vx = \bar{\vx}$,
the remaining {\em subproblem} becomes much easier to solve, perhaps by decoupling into smaller problems.  In our study, the smaller problems will correspond to scenarios and facilities.
The subproblem has the form
\begin{equation}
\label{LBBD:subproblem}
\mathrm{SP}(\bar{\vx}) = \min_{\vy} \big\{ f(\bar{\vx},\vy) \;\big|\; C(\bar{\vx},\vy), \; \vy\in D_y \big\}
\end{equation}
The key to LBBD is analyzing the subproblem solution so as to find a function $B_{\bar{\vx}}(\vx)$ that provides a lower bound on $f(\vx,\vy)$ for any given $\vx\in D_{x}$.  The bound must be sharp for $\vx=\bar{\vx}$; that is, $B_{\bar{\vx}}(\bar{\vx})=\mathrm{SP}(\bar{\vx})$.  The bounding function is derived from the {\em inference dual} of the subproblem in a manner discussed below.  In classical Benders decomposition, the subproblem is an LP problem, and the inference dual is the LP dual.  

Each iteration of the LBBD algorithm begins by solving a {\em master problem}:
\begin{equation}
\label{eq:relaxMaster}
\mathrm{MP}(\bar{\vX}) = \min_{\vx,\beta} \big\{ \beta \;\big|\; \beta\geq B_{\bar{\vx}}(\vx), \; \mbox{all}\; \bar{\vx}\in \bar{\vX}; \;\; \vx\in D_x \big\}
\end{equation}
where the inequalities $\beta\geq B_{\bar{\vx}}(\vx)$ are {\em Benders cuts} obtained from previous solutions $\bar{\vx}$ of the subproblem.  There may be several cuts for a given $\bar{\vx}$, but for simplicity we assume in this section there is only one.  Initially, the set $\bar{\vX}$ can be empty, or it can contain a few solutions obtained heuristically to implement a ``warm start.''  
The optimal value MP$(\bar{\vX})$ of the master problem is a lower bound on the optimal value of the original problem (\ref{LBBD:general}).  If $\bar{\vx}$ is an optimal solution of the master problem, the corresponding subproblem is then solved to obtain SP$(\bar{\vx})$, which is an upper bound on the optimal value of (\ref{LBBD:general}).  A new Benders cut $\beta\geq B_{\bar{\vx}}(\vx)$ is generated for the master problem and $\bar{\vx}$ added to $\bar{\vX}$ in (\ref{eq:relaxMaster}).  The process repeats until the lower and upper bounds provided by the master problem and subproblem converge; that is, until
MP$(\bar{\vX})=\min_{\bar{\vx}\in\bar{\vX}} \{\mathrm{SP}(\bar{\vx})\}$.  The following is proved in \cite{Hoo00}:
\begin{theorem}
If $D_x$ is finite, the LBBD algorithm converges to an optimal solution of (\ref{LBBD:general}) after a finite number of iterations.
\end{theorem}

The inference dual of the subproblem seeks the tightest bound on the objective function that can be inferred from the constraints.  Thus the inference dual is
\begin{equation}
\mathrm{DSP}(\bar{\vx}) = \max_{P\in\mathcal{P}} \Big\{ \gamma \; \Big| \; \big(C(\bar{\vx}),\; y\in D_y\big) \stackrel{P}{\Rightarrow} \big( f(\bar{\vx},\vy)\geq\gamma)\big) \Big\}
\label{eq:dual}
\end{equation}
where $A\stackrel{P}{\Rightarrow}B$ indicates that proof $P$ deduces $B$ from $A$.  The inference dual is always defined with respect to set $\mathcal{P}$ of valid proofs.  In classical linear programming duality, valid proofs consist of nonnegative linear combinations of the inequality constraints in the problem.  We assume a strong dual, meaning that SP$(\bar{\vx})=\mathrm{DSP}(\bar{\vx})$.  The dual is strong when the inference method is complete.  For example, the classical Farkas Lemma implies that nonnegative linear combination is a complete inference method for linear inequalities.  Indeed, any exact optimization method is associated with a complete inference method that it uses to prove optimality, perhaps one that involves branching, cutting planes, constraint propagation, and so forth.  

In the context of LBBD, the proof $P$ that solves the dual (\ref{eq:dual}) is the proof of optimality the solver obtains for the subproblem (\ref{LBBD:subproblem}).  The bounding function $B_{\bar{\vx}}(\vx)$ is derived by observing what bound on the optimal value {\em this same proof} $P$ can logically deduce for a given $\vx$, whence the description ``logic-based.'' In practice, the solver may not reveal how it proved optimality, or the proof may be too complicated to build a useful cut.  One option in such cases is to tease out the structure of the proof by re-solving the subproblem for several values of $\vx$ and observing the optimal value that results.  This information can be used to design {\em strengthened nogood cuts} that provide useful bounds for many values of $\vx$ other than $\bar{\vx}$.  Another approach is to use {\em analytical Benders cuts}, which deduce bounds on the optimal value when $\bar{\vx}$ is changed in certain ways, based on structural characteristics of the subproblem and its current solution.  We will employ both of these options.

{\em Branch and check} is a variation of LBBD that solves the master problem only once and generates Benders cuts on the fly.  It is most naturally applied when the master problem is solved by branching. Whenever the branching process discovers a solution $\bar{\vx}$ that is feasible in the current master problem, the corresponding subproblem is solved to obtain one or more Benders cuts, which are added to the master problem.  Branching then continues and terminates in the normal fashion, all the while satisfying Benders cuts as they accumulate.  Branch and check can be superior to standard LBBD when the master problem is much harder to solve than the subproblems.  

A common enhancement of LBBD and other Benders methods is a {\em warm start}, which includes initial Benders cuts in the master problem.  Recent studies that benefit from this technique include \cite{angulo2016improving}, \cite{elcci2018chance}, and \cite{HecHooKim19}.  Benders cuts can also be aggregated before being added to the master problem, a technique first explored in \cite{birge1988multicut}.  A particularly useful enhancement for LBBD is to include a relaxation of the subproblem in the master problem, where the relaxation is written in terms of the master problem variables \citep{
hooker2007planning,fazel2012using}.  We employ this technique in the present study.

\section{Benders Formulation of Planning \& Scheduling}

We apply LBBD to the generic planning and scheduling problem by placing the assignment decision in the master problem and the scheduling decision in the subproblem.  The master problem is therefore
\[
\min_{\vx} \Big\{ g(\vx) + \sum_{\omega\in\Omega} \pi_{\omega}\beta_{\omega} \; \Big| \; \mbox{Benders cuts}; \; \mbox{subproblem relaxation}; \; x_j\in I, \;\mbox{all}\; j\in J \Big\}
\]
The Benders cuts provide lower bounds on each $\beta_{\omega}$.  The cuts and subproblem relaxation are somewhat different for each variant of the problem we consider below.  The scheduling subproblem decouples into a separate problem for each facility and scenario.  If $\bar{\vx}$ is an optimal solution of the master problem, the scheduling problem for facility $i$ and scenario $\omega$ is
\[
\mathrm{SP}_{i\omega}(\bar{\vx}) = 
\min_{\bm{s}} \Big\{ h_i(\bm{s},\bar{\vx},\omega) \;\Big| \; s_j\in [r_j, d_j-p^{\omega}_{ij}], \; \mbox{all}\; j\in J_i(\bar{\vx}); 
\hspace{-2ex} \sum_{\substack{j\in J_i(\bar{\vx})\\0\leq t\leq s_j+p^{\omega}_{ij}}} \hspace{-3ex} c_{ij} \leq K_i, \;\mbox{all}\; t \Big\}
\]

We solve the master problem and subproblem by formulating the former as an MILP problem and the latter as a CP problem.  
In the master problem, we let variable $x_{ij}=1$ when task $j$  is assigned to facility $i$.  The master problem becomes
\begin{equation}
\begin{array}{ll}
\mbox{minimize} & 
{\ds
\hat{g}(\vx) + \sum_{\omega\in\Omega} \pi_{\omega}\beta_{\omega}
} \vsp \\
\mbox{subject to} & 
{\ds
\sum_{i\in I} x_{ij} = 1, \;\; j\in J
} \vsp \\
& \mbox{Benders cuts} \vsp \\
& \mbox{subproblem relaxation} \vsp \\
& x_{ij}\in \{0,1\}, \; i\in I, \; j\in J
\end{array} \label{eq:MILP}
\end{equation}
where $\vx$ now denotes the matrix of variables $x_{ij}$.  If $\bar{\vx}$ is an optimal solution of the master problem, the subproblem for each facility $i$ and scenario $\omega$ becomes
\begin{equation}
\begin{array}{ll}
\text{minimize}  & \hat{h}_i(\bm{s},\bar{\vx},\omega) \vsp \\
\text{subject to} & \text{cumulative}\Big(\big(s_j\,\big|\,j \in J_i(\bar{\vx})\big), \,
\big(p^{\omega}_{ij}\,\big|\,j \in J_i(\bar{\vx})\big), \,\big(c_{ij}\,\big|\,j \in J_i(\bar{\vx})\big), \, 
K_i \Big) \vsp \\
& s_j \in [r_j, d_i-p^{\omega}_{ij} ], \; j \in J_i(\bar{\vx})
\end{array} \label{eq:CP}
\end{equation}
The optimal value of (\ref{eq:CP}) is again SP$_{i\omega}(\bar{\vx})$.  The cumulative global constraint in (\ref{eq:CP}) is a standard feature of CP models and requires that the total resource consumption at any time on facility $i$ be at most $K_i$.  

To solve a problem (\ref{eq:P&Srisk}) that incorporates risk, one need only replace the objective function of (\ref{eq:MILP}) with 
\[
\hat{g}(\vx) + (1-\lambda)\sum_{\omega\in\Omega} \pi_{\omega}\beta_{\omega} + \lambda\mathbb{D}_{\omega} [\beta_{\omega}]
\]
and otherwise proceed as in the risk-neutral case.

\subsection{Minimum Makespan Problem}
\label{Sec:Algorithm:Makespan}

We begin by considering a minimum makespan problem in which the tasks have release times and no deadlines.  
The first-stage objective function is $g(\vx)=0$, and so we have $\hat{g}(\vx)=0$ in the MILP model (\ref{eq:MILP}).  The second-stage objective function is the finish time of the last task to finish:
\[
h(\bm{s},\vx,\omega) = \max_{j\in J} \Big\{ s_j + p^{\omega}_{x_j j} \Big\}
\]
This objective function is incorporated into the CP problem (\ref{eq:CP}) by setting $\hat{h}_i(\bm{s},\bar{\vx},\omega)=M$ and adding to (\ref{eq:CP}) the constraints $M\geq s_j + p^{\omega}_{ij}$ for all $j\in J_i(\bar{\vx})$.  Since there are no deadlines, we assume $d_j=\infty$ for all $j\in J$.

Both strengthened nogood cuts and analytic Benders cuts can be developed for this problem.  A simple nogood cut for scenario $\omega$ can take the form of a set of inequalities
\begin{equation}
\beta_{\omega} \geq \beta_{i\omega}, \; i\in I
\label{eq:nogood0}
\end{equation}
where each $\beta_{i\omega}$ is bounded by
\begin{equation}
\beta_{i\omega} \geq \mathrm{SP}_{i\omega}(\bar{\vx}) \Big( \sum_{j\in J_i(\bar{\vx})} \hspace{-1.5ex}  x_{ij} - |J_i(\bar{\vx})| + 1 \Big)
\label{eq:nogood1}
\end{equation}
and where $\bar{\vx}$ is the solution of the current master problem.  The cut says that if all the jobs in $J_i(\bar{\vx})$ are assigned to facility $i$, possibly among other jobs, then the makespan of facility $i$ in scenario $\omega$ is at least the current makespan SP$_{i\omega}(\bar{\vx})$.  The cut is weak, however, because if even one job in $J_i(\bar{\vx})$ is not assigned to $i$, the bound in (\ref{eq:nogood1}) becomes useless.  The cut can be strengthened by heuristically assigning proper subsets of the jobs in $J_i(\bar{\vx})$ to facility $i$, and re-computing the minimum makespan for each subset, to discover a smaller set of jobs that yields the same makespan.  This partially reveals which task assignments serve as premises of the optimality proof.  Then $J_i(\bar{\vx})$ in (\ref{eq:nogood1}) is replaced with this smaller set to strengthen the cut.  This simple scheme, and variations of it, can be effective when the makespan problem solves quickly \cite{hooker2007planning}.

A stronger cut can be obtained without re-solving the makespan problem by using an analytical Benders cut.  We introduce a cut based on the following lemma:

\begin{lemma} \label{lemma1}
Consider a minimum makespan problem $P$
in which each task $j\in J$ has release time $r_j$ and processing time $p_j$, with no deadlines.  Let $M^*$ denote the minimum makespan for $P$,  and $\hat{M}$ the minimum makespan for the problem 
$\hat{P}$ that is identical to $P$ except that the tasks
in a nonempty set $\hat{J}\subset J$ are removed. Then
\begin{equation}
M^* - \hat{M} \leq \Delta + r^+ - r^-
\label{lemma_eq}
\end{equation}
where $\Delta = \sum_{j \in \hat{J}} p_j$, $r^+ = \max_{j \in J} \{ r_j\}$ 
is the latest release time, and $r^- = \min_{j \in J} \{ r_j\}$ is the earliest release time.
\end{lemma}

\begin{proof}
Consider any solution of $\hat{P}$ with makespan $\hat{M}$. We will construct a feasible solution for $P$ by extending this solution.
If $\hat{M} > r^+$, we schedule all the tasks in $\hat{J}$ sequentially starting from time $\hat{M}$, resulting in makespan $\hat{M}+\Delta$.
This is a feasible solution for $P$, and we have $M^* \leq \hat{M} + \Delta$.  The lemma follows because $r^+ - r^-$ is nonnegative.
If $\hat{M} < r^+$, we schedule all the tasks in $\hat{J}$ sequentially starting from time $r^+$ to obtain a solution with makespan of $r^+ + \Delta$. Again this is a feasible solution for $P$, and we have
$M^* \leq r^+ + \Delta$. This implies
\begin{equation}
M^* - \hat{M} \leq r^+ - \hat{M} + \Delta
\label{proof_eq}
\end{equation}
Because $\hat{M}$ is at least $r^-$, \eqref{proof_eq} implies \eqref{lemma_eq}, and the lemma follows. \qed
\end{proof}


We can now derive a valid analytical cut:
\begin{theorem}
Inequalities (\ref{eq:nogood0}) and the following comprise a valid Benders cut for scenario $\omega$:
\begin{equation}
\label{eq:makespanCutA}
\beta_{i\omega} \geq
\left\{
\begin{array}{ll}
{\ds
	\mathrm{SP}_{i\omega}(\bar{\vx}) - \Big( \hspace{-0.5ex} \sum_{j \in J_i(\bar{\vx})} \hspace{-1ex} (1 - x_{ij}) p^{\omega}_{ij} + r^+ - r^- \Big),
} & \text{if $x_{ij}=0$ for some $j\in J_i(\bar{\vx})$} \vsp \\
\mathrm{SP}_{i\omega} (\bar{\vx}),  & \text{otherwise}
\end{array} 
\right\}, \; i\in I 
\end{equation}
\end{theorem}

\begin{proof}
The cut clearly provides a sharp bound $\max_{i\in I} \{\mathrm{SP}_{i\omega} (\bar{\vx})\}$ when $\vx=\bar{\vx}$, because the second line of (\ref{eq:makespanCutA}) applies in this case.  The validity of the cut follows immediately from Lemma~\ref{lemma1}. \qed
\end{proof}

We linearize the cut (\ref{eq:makespanCutA}) as follows:
\begin{equation} 
\label{Cut:Makespan:Strong}
\begin{array}{l}
{\ds
\beta_{i\omega} \geq \mathrm{SP}_{i\omega}(\bar{\vx}) - \hspace{-1ex} \sum_{j \in J_i(\bar{\vx})} \hspace{-1ex} (1 - x_{ij}) p^{\omega}_{ij} - z_{i\omega}
} \vsp \\
{\ds
z_{i\omega} \leq \left( r^+ - r^- \right) \hspace{-1ex} \sum_{j \in J_i(\bar{\vx})} \hspace{-1ex} (1 - x_{ij})
} \vsp \\
z_{i\omega} \leq r^+ - r^-
\end{array} 
\end{equation}
The Benders cut is inserted into the master problem by including inequalities (\ref{Cut:Makespan:Strong}) for each $i\in I$ and $\omega\in\Omega$, along with the inequalities (\ref{eq:nogood0}).

The inequalities (\ref{Cut:Makespan:Strong}) incur the expense of introducing a new continuous variable $z_{i\omega}$ for each $i$ and $\omega$, which may not be 
desirable while solving large-scale problems. As an alternative,
a slightly weaker Benders cut can be used.
\begin{corollary}
Inequalities (\ref{eq:nogood0}) and the following comprise a valid Benders cut for scenario $\omega$:
\begin{equation} 
\label{Cut:Makespan:Strong:ver2}
\beta_{i\omega}  \geq \mathrm{SP}_{i\omega}(\bar{\vx}) - \hspace{-1ex} \sum_{j \in J_i(\bar{\vx})} \hspace{-1ex} (1 - x_{ij}) p^{\omega}_{ij}
- \left( r^+ - r^- \right) |J_i(\bar{\vx})|^{-1} \hspace{-1ex} \sum_{j \in J_i(\bar{\vx})} \hspace{-1ex} (1 - x_{ij}), \;\; i\in I
\end{equation}
\end{corollary}

\begin{proof}
We first note that the inequalities provide a sharp bound if $\vx=\bar{\vx}$, because in this case $x_{ij}=1$ for all $j\in J_i(\bar{\vx})$, and (\ref{Cut:Makespan:Strong:ver2}) is identical to the second line of (\ref{Cut:Makespan:Strong}).  If $x_{ij}=0$ for some $j\in J_i(\bar{\vx})$, we have 
\[
\left( r^+ - r^- \right) |J_i(\bar{\vx})|^{-1} \hspace{-1ex} \sum_{j \in J_i(\bar{\vx})} \hspace{-1ex} (1 - x_{ij}) \leq  r^+ - r^-
\]
because $\sum_{j\in J_i(\bar{\vx})} (1-x_{ij})\leq J_i(\bar{\vx})$.  Thus (\ref{Cut:Makespan:Strong:ver2}) is implied by the first line of (\ref{Cut:Makespan:Strong}) and is therefore valid. \qed
\end{proof}

Finally, we add a subproblem relaxation to the master problem.  We use a relaxation from \cite{hooker2007planning},
modified to be scenario-specific:
\begin{equation}
\beta_{i\omega} \geq \frac{1}{K_i} \sum_{j \in J} c_{ij} p^{\omega}_{ij} x_{ij},
\quad i \in I, \; \omega\in \Omega
\label{eq:relaxMakespan}
\end{equation}
This relaxation is valid for arbitrary release times and deadlines.

\subsection{Minimum cost problem}
\label{Sec:Algorithm:MinCost} 
In the minimum cost problem, there is only a fixed cost $\phi_{ij}$ associated with assigning task $j$ to facility $i$.  So we have 
\[
\hat{g}(\vx) = \sum_{i\in I} \sum_{j\in J} \phi_{ij}x_{ij}
\]
in the MILP master problem (\ref{eq:MILP}), and we set $\beta_{\omega}=0$ for $\omega\in\Omega$.  The subproblem decouples into a feasibility problem for each $i$ and $\omega$, because $\hat{h}_i(\bm{s},\bar{\vx},\omega) = 0$.  

A Benders cut is generated for each $i$ and $\omega$ when the corresponding scheduling problem (\ref{eq:CP}) is infeasible.  A simple nogood cut is
\begin{equation}
\sum_{j\in J_i(\bar{\vx})} (1-x_{ij}) \geq 1
\end{equation}
The cut can be strengthened in a manner similar to that used for the makespan problem.

To create a subproblem relaxation for the master problem, one can  exploit the fact that we now have two-sided time windows $[r_j,d_j]$.  Let $J(t_1,t_2)$ be the set of tasks $j$ for which $[r_j,d_j]\subseteq [t_1,t_2]$.  Adapting an approach from \citep{hooker2007planning}, one can add the following inequalities to the master problem for each $i\in I$:
\begin{equation}
\label{Relaxation:Cost}
\frac{1}{K_i} \sum_{j \in J(t_1,t_2)} 
\hspace{-2ex} p^{\min}_{ij} c_{ij} x_{ij} \leq t_2 - t_1, \;\;
t_1\in \{\bar{r}_1,\ldots,\bar{r}_{n'}\}, \; t_2\in \{\bar{d}_1,\ldots,\bar{d}_{n''}\}
\end{equation}
where $\bar{r}_1,\ldots,\bar{r}_{n'}$ are the distinct release times among $r_1,\ldots, r_n$, and $\bar{d}_1,\ldots,\bar{d}_{n''}$ the distinct deadlines among $d_1,\ldots,d_n$.  Some of these inequalities may be redundant, and a method for detecting them is presented in \citep{hooker2007planning}.  Because the relaxation must be valid across all scenarios, the processing time is set to $p^{\min}_{ij}=\min_{\omega\in\Omega} \{p^{\omega}_{ij}\}$.

\subsection{Minimum tardiness problem}
\label{Sec:Algorithm:Tardiness}
In this section, we consider a minimum tardiness problem in which tasks are all released at time zero but have different due dates $\bar{d}_j$. There are no hard deadlines, and so we let $d_j=\infty$ for all $j\in J$. As in the minimum makespan problem, there is no first-stage cost, so that $\hat{g}(\vx)=0$ in the MILP model (\ref{eq:MILP}).  The second-stage objective function is expected total tardiness, and we have 
\[
\hat{h}_i(\bm{s},\vx,\omega) = \hspace{-1ex} \sum_{j\in J_i(\vx)} \hspace{-1ex} \big( s_j+p^{\omega}_{ij} - \bar{d}_j \big)^+ 
\]
in the CP scheduling problem (\ref{eq:CP}).  Here $\alpha^+ = \max\{0,\alpha\}$.

The following analytic Benders cut can be adapted from \cite{Hoo12}:
\[
\beta_{\omega} \geq
\sum_{i\in I} \Big( \mathrm{SP}_{i\omega}(\bar{\vx}) -
\hspace{-1ex} \sum_{j\in J_i(\bar{\vx})} \hspace{-1ex} \Big( \hspace{-0.5ex} \sum_{j'\in J_i(\bar{\vx})} \hspace{-1.5ex} p^{\omega}_{ij'} - \bar{d}_j\Big)^+ (1-x_{ij}) \Big)
\]
The cut is added to (\ref{eq:MILP}) for each $\omega\in\Omega$.  Strengthened nogood cuts similar to those developed for the makespan problem can also be used.

Two subproblem relaxations can be adapted from \cite{hooker2007planning}.  The simpler one is analogous to (\ref{Relaxation:Cost}) and adds the following inequalities to (\ref{eq:MILP}) for each $i$ and $\omega$
\[
\beta_{i\omega} \geq \frac{1}{K_i} \hspace{-0.5ex} \sum_{j'\in J(0,\bar{d}_j)} \hspace{-2ex} p^{\min}_{ij'}c_{ij'}x_{ij'} - \bar{d}_j, \;\; j\in J
\]
along with the bounds $\beta_{i\omega}\geq 0$.  A second relaxation more deeply exploits the structure of the subproblem.  For each facility $i$ and scenario $\omega$, let $\tau^{\omega}_i$ be a permutation of $\{1,\ldots,n\}$ such that $p^{\omega}_{i\tau^{\omega}_i(1)}c_{i\tau^{\omega}_i(1)} \leq \cdots \leq p^{\omega}_{i\tau^{\omega}_i(n)}c_{i\tau^{\omega}_i(n)}$.  We also assume that tasks are indexed so that $\bar{d}_1\leq \cdots\leq \bar{d}_n$.  Then we add the following inequalities to the master problem (\ref{eq:MILP}) for each $i$ and $\omega$:
\[
\beta_{i\omega} \geq
\frac{1}{K_i} \sum_{j'\in J} p^{\omega}_{i\tau^{\omega}_i(j')}c_{i\tau^{\omega}_i(j')}x_{i\tau^{\omega}_i(j')} - \bar{d}_j - (1-x_{ij})U_{ij\omega}, \;\; j\in J
\]
where
\[
U_{ij\omega} = \frac{1}{K_i} \sum_{j'\in J} p^{\omega}_{i\tau^{\omega}_i(j')}c_{i\tau^{\omega}_i(j')} - \bar{d}_j
\]

\section{The Integer L-Shaped Method}
\label{Sec:intLShaped}
The integer L-Shaped method is a Benders-based algorithm proposed by \cite{laporte1993integer} to solve two-stage stochastic integer programs. It terminates in finitely many iterations when the problem has complete recourse and binary first-stage variables. It is similar to branch and check in that Benders cuts are generated while solving the first-stage problem by branching.  It differs in that it uses subgradient cuts derived from a linear programming relaxation of the subproblem rather than combinatorial cuts derived from the original subproblem.  It also uses a simple integer nogood cut to ensure convergence, but the cut is quite weak and does not exploit the structure of the subproblem as does branch and check.  We describe the integer \mbox{L-shaped} method here as it applies to minimizing makespan in the planning and scheduling problem.

We first state an MILP model of the deterministic equivalent problem, as it will play a benchmarking role in computational testing.  We index discrete times by $t\in T$ and introduce a 0--1 variable $z^{\omega}_{ijt}$ that is 1 if task $j$ starts at time $t$ on facility $i$ in scenario $\omega$.  The model is
\begin{equation}
\begin{array}{lll}
\text{minimize} & 
{\ds 
\sum_{\omega \in \Omega} \pi_{\omega} \beta_{\omega} 
} & (a)  \vsp \\
\text{subject to} &
{\ds
\sum_{i\in I} x_{ij} = 1, \;\; j\in J
} & (b) \vsp \\
& \beta_{\omega} \geq \beta_{i\omega}, \;\; i\in I, \; \omega\in\Omega & (c) \vsp \\
& x_{ij}\in \{0,1\}, \;\; i\in I, \; j\in J & (d) \vsp \\
& {\ds
\beta_{\omega} \geq \sum_{t\in T} (t + p^{\omega}_{ij}) z^{\omega}_{ijt}, \;\; i\in I, \; j\in J, \; \omega\in\Omega 
} & (e) \vsp \\
& z^{\omega}_{ijt} \leq x_{ij}, \;\; i\in I, \; j\in J, \; t\in T, \; \omega\in\Omega & (f) \vsp \\
& {\ds 
\sum_{i\in I} \sum_{t\in T} z^{\omega}_{ijt} = 1, \;\; j\in J, \;  \omega\in\Omega 
} & (g) \vsp \\
& {\ds
\sum_{j\in J} \sum_{t'\in T^{\omega}_{tij}} \hspace{-1ex} c_{ij} z^{\omega}_{ijt'} \leq K_i, \;\; i\in I, \; t\in T, \; \omega\in\Omega
} & (h) \vsp \\
& z^{\omega}_{ijt}=0, \;\; i\in I, \; \omega\in \Omega, \; j\in J, \; \mbox{all}\; t\in T\; \mbox{with} \; t<r_j & (i) \vsp \\
& z^{\omega}_{ijt} \in \{0,1\}, \;\; i\in i, \; j\in J, \; t\in T, \; \omega\in\Omega & (j)
\end{array} \label{eq:DEQ}
\end{equation}
where $T^{\omega}_{tij}= \{t'\;|\; 0\leq t' \leq t-p^{\omega}_{ij}\}$.  In the integer L-shaped method, the first-stage minimizes (\ref{eq:DEQ}a) subject to (\ref{eq:DEQ}b)--(\ref{eq:DEQ}d)
and Benders cuts that provide bounds on $\beta_{i\omega}$.  The Benders cuts consist of classical Benders cuts derived from the linear relaxation of the second-stage scheduling problem for each $i$ and $\omega$, as well as integer cuts.  If $\bar{\vx}$ is an optimal solution of the first-stage problems, the second-stage problem for facility $i$ and scenario $\omega$ is
\begin{equation}
\begin{array}{ll}
\text{minimize} & M \vsp \\
\text{subject to} &
{\ds
M \geq \sum_{t\in T} (t + p^{\omega}_{ij}) z^{\omega}_{ijt}, \;\; j\in J_i(\bar{\vx}) 
} \vsp \\
& {\ds 
\sum_{t\in T} z^{\omega}_{ijt} = 1, \;\; j\in J_i(\bar{\vx}) 
} \vsp \\
& {\ds
\sum_{j\in J} \sum_{t'\in T^{\omega}_{tij}} \hspace{-1ex} c_{ij} z^{\omega}_{ijt'} \leq K_i, \;\; t\in T
}  \vsp \\
& x^{\omega}_{ijt} \in \{0,1\}, \;\; j\in J_i(\bar{\vx}), \; t\in T \vsp \\
& z^{\omega}_{ijt}=0, \;\; j\in J_i(\bar{\vx}), \; \mbox{all}\; t\in T\; \mbox{with} \; t<r_j  \vsp \\
& z^{\omega}_{ijt} \in \{0,1\}, \;\; j\in J_i(\bar{\vx}), \; t\in T
\end{array} \label{eq:LshapedSub}
\end{equation}
The following integer cut is used for each $i$ and $\omega$ to ensure convergence:
\begin{equation}
\beta_{i\omega} \geq \big( \mathrm{SP}_{i\omega}(\bar{\vx}) - \mathrm{LB}\big) \Big( \hspace{-0.5ex} \sum_{j\in J_i(\bar{\vx})}\hspace{-1.5ex} x_{ij} - \hspace{-1ex} \sum_{j\not\in J_i(\bar{\vx})} \hspace{-1.5ex} x_{ij} - |J_i(\bar{\vx})| + 1 \Big) + \mathrm{LB}
\label{eq:integerCut}
\end{equation}
where LB is a global lower bound on makespan.  Note that if $\mathrm{LB}=0$, (\ref{eq:integerCut}) is weaker than the unstrengthened nogood cut (\ref{eq:nogood1}).  This is because (\ref{eq:integerCut}) becomes useless if $x_{ij}\neq \bar{x}_{ij}$ for even one $j\in J$, while (\ref{eq:nogood1}) becomes useless only if $x_{ij}\neq \bar{x}_{ij}$ for some $j\in J_i(\bar{\vx})$.  If a bound  $\mathrm{LB}>0$ is available, (\ref{eq:integerCut}) is still weaker than (\ref{eq:nogood1}) if the same bound is added to (\ref{eq:nogood1}) by writing
\[
\beta_{i\omega} \geq \big( \mathrm{SP}_{i\omega}(\bar{\vx}) - \mathrm{LB}\big) \Big( \sum_{j\in J_i(\bar{\vx})} \hspace{-1.5ex}  x_{ij} - |J_i(\bar{\vx})| + 1 \Big) + \mathrm{LB}
\]

Finally, we strenghthen the initial master problem by adding bounds of the form
\begin{equation}
\beta_{\omega} \geq \beta^{\mathrm{LB}}_{\omega}, \;\; \omega\in\Omega.
\label{eq:initialBounds}
\end{equation}
Here $\beta^{\mathrm{LB}}_{\omega}$ is a lower bound on makespan obtained by solving the LP relaxation of (\ref{eq:DEQ}) for fixed scenario $\omega$.  We use the same bound in LBBD and branch-and-check methods.

\section{Computational Study}
\label{Sec:CompStudy}
%
In this section, we describe computational experiments we conducted for the stochastic planning and scheduling
problem, with the objective of minimizing makespan.  All experiments are conducted on a personal computer with a 2.80 GHz 
Intel\textsuperscript{\textregistered} Core\textsuperscript{\texttrademark} i7-7600 processor
and 24 GB memory running on a Microsoft Windows 10 Pro.
All MILP and CP formulations are solved in \texttt{C++} using the \texttt{CPLEX} and \texttt{CP Optimizer} engines of
\texttt{IBM}\textsuperscript{\textregistered} \texttt{ILOG}\textsuperscript{\textregistered}   \texttt{CPLEX}\textsuperscript{\textregistered} \texttt{12.7 Optimization Studio},
respectively.  We use only use a single thread in all computational experiments. We modify \texttt{CP Optimizer} parameters to execute an extended filtering and DFS search. The rest of the parameters are set to their default values for both \texttt{CPLEX} and \texttt{CP Optimizer} engines. Lastly, we use the \texttt{Lazy Constraint Callback} function of \texttt{CPLEX} to implement branch and check.

\subsection{Problem Instances}
\label{Sec:CompStudy:InstanceGen}

We generate problem instances by combining ideas from \cite{hooker2007planning} and \cite{atakan2017minimizing}.
We first generate the deterministic problem as in \cite{hooker2007planning}. Let $|I| = m$ and $|J| = n$.
The capacity limits of the facilities is set to $K_i = 10$ for all $i\in I$, and integer capacity requirements of tasks are drawn from a uniform distribution on $[1,10]$.  Integer release times are drawn from a uniform distribution on $[0,~2.5 n(m+1)/m]$.  For each facility $i\in I$, integer mean processing times $\bar{p}_{ij}$ are drawn from a uniform distribution on $[2, ~25 - 10(i-1)/(m-1)]$. This causes facilities with a higher index to process tasks more rapidly.  

We then follow \cite{atakan2017minimizing} by perturbing the mean processing times to obtain a set of scenarios. In particular, we first divide the tasks into two groups, one group containing tasks $i$ for which $0<\bar{p}_{ij} \leq 16$, and the other group containing the remainder of the tasks.  We then generate a perturbation parameter $\epsilon^{\omega}$ for each scenario $\omega \in \Omega$ from a mixture of uniform distributions.  Specifically, for tasks in the first group, $\epsilon^{\omega}$ is distributed uniformly on the interval $[-0.1,~ 0.5]$ with probability 0.9 and on the interval $[2.0, ~3.0]$ with probability 0.1.  For tasks in the second group, $\epsilon^{\omega}$ is distributed uniformly on the
interval $[-0.1,~ 0.5]$ with probability 0.99 and on the interval $[1.0, ~1.5]$ with probability 0.01.
Finally, we generate the processing times under scenario $\omega \in \Omega$ by letting
$p^{\omega}_{ij}= \lceil \bar{p}_{ij}(1+\epsilon^{\omega}) \rceil$.

\subsection{Computational performance}
\label{Sec:CompStudy:Performance}
In this section, we report comparisons of LBBD and branch and check with the integer \mbox{L-shaped} method.  The experiments are designed to investigate how various algorithmic features affect performance.  All results presented in the tables to follow are averages over 3 random instances. 

Table~\ref{table:CompPerf1} compares computation times and optimality gaps for seven solution methods.  Each method solves the first-stage problem using the MILP engine in \texttt{CPLEX}.   \begin{itemize} 
	
	\item {\em Deterministic equivalent MILP.}  We solve the deterministic equivalent model (\ref{eq:DEQ}) using the MILP engine in \texttt{CPLEX}.
	
	\item {\em Standard integer L-shaped method.} We decouple the second-stage problem by facility and scenario and solve the resulting problems and their LP relaxations using the MILP engine of \texttt{CPLEX} whenever a candidate incumbent solution is identified. We then we add the integer cut \eqref{eq:integerCut} and the classical Benders cut from the LP relaxation for each facility and scenario.  The initial bounds (\ref{eq:initialBounds}) are included in the master problem, even though they are not standard, because previous experience indicates that they significantly enhance performance.  The subproblem relaxation (\ref{eq:relaxMakespan}) is likewise included in the master problem for fair comparison with LBBD and branch and check, where it is standard.
	
	\item {\em Integer L-shaped method with CP.}  We modify the standard method by solving the second-stage subproblems with CP rather than MILP. Integer cuts are as before, and classical Benders cuts are derived from the LP relaxation of the MILP model as before.  The initial bounds (\ref{eq:initialBounds}) and subproblem relaxation (\ref{eq:relaxMakespan}) are again included in the master problem.
	
	\item {\em Standard LBBD with nogood cuts.}  We use (\ref{eq:nogood0}) and unstrengthened nogood cuts (\ref{eq:nogood1}). We solve the decoupled subproblems by \texttt{CP Optimizer}.  The initial bounds (\ref{eq:initialBounds}) are included in the master problem for comparability with the integer \mbox{L-shaped} method.
	
	\item {\em Standard LBBD with analytical cuts.}  We use (\ref{eq:nogood0}) and  analytical cuts (\ref{Cut:Makespan:Strong:ver2}) rather than nogood cuts.  The decoupled subproblems are solved by \texttt{CP Optimizer}.  The initial bounds (\ref{eq:initialBounds}) are again included in the master problem.
	 
	\item {\em Branch and check with nogood cuts.} We use (\ref{eq:nogood0}) and unstrengthened nogood cuts (\ref{eq:nogood1}). We solve the decoupled subproblems by \texttt{CP Optimizer}.  The initial bounds (\ref{eq:initialBounds}) are included in the master problem.
	
	\item {\em Branch and check with analytical cuts.} We use (\ref{eq:nogood0}) and  analytical cuts (\ref{Cut:Makespan:Strong:ver2}) rather than nogood cuts. The decoupled subproblems are solved by \texttt{CP Optimizer}.  The initial bounds (\ref{eq:initialBounds}) are again included in the master problem.
	
\end{itemize}


In addition to average computation time (in seconds), Table~\ref{table:CompPerf1} reports the optimality gap obtained for each solution method, defined as $(\mathrm{UB}-\mathrm{LB})/\mathrm{UB}$.  For the deterministic equivalent and branch-and-check methods, UB and LB are, respectively, the upper and lower bounds obtained from \texttt{CPLEX} upon solution of the master problem.  For standard LBBD, UB and LB are, respectively, the smallest subproblem optimal value and the largest master problem optimal value obtained during the Benders algorithhm.

\begin{table}[!t]
\small
\centering
\caption{Average computation time in seconds over 3 instances (upper half of table) and average relative optimality gap (lower half) for various solution methods, based on 10 tasks and 2 facilities.}
\label{table:CompPerf1}
\vspace{3ex}
\begin{tabular}{r|r@{}r@{}r@{}r@{}r@{}r@{}rrrr}
	&Determ.&&Integer&&Integer&& LBBD  & LBBD  & B\&Ch & B\&Ch \\
	Scenarios	
	&equiv.	&&L-shaped	
	&&L-shaped
	&&Nogood &Analytic
	&Nogood	&Analytic \\
	& MILP	&&method	&&with CP&&cuts	&cuts	&cuts	&cuts \\
	\cmidrule{1-11} 
1&  2.4 & &  127.3&&  27.9&&  2.0&  0.6&  1.8&  1.6 \\
5&  475.8 &$^{\dagger\dagger}$&  839.2&&  149.3&&  12.1&  3.0&  3.3&  1.7 \\
10&  *& & 2316.9&$^{\dagger}$&  437.8&&  27.4&  7.3&  5.1&  2.8 \\
50&  *& &*&&  2517.8&$^{\dagger\dagger}$&  243.1&  42.8&  33.3&  18.0 \\
100&  *&& *&& *&&  952.8&  118.8&  80.1&  30.2 \\
500&  *&& *&& *&&  *&  900.9&  416.2&  166.2 \\
	\cmidrule{1-11} 
1&  0.0&&  0.0&&  0.0&&  0.0&  0.0&  0.0&  0.0 \\
5&  7.8&&  0.0&&  0.0&&  0.0&  0.0&  0.0&  0.0 \\ 
10&  12.4&&  3.8&&  0.0&&  0.0&  0.0&  0.0&  0.0 \\
50&  17.4&&  21.7&&  13.9&&  0.0&  0.0&  0.0&  0.0 \\
100&  25.4&&  25.4&&  21.7&&  0.0&  0.0&  0.0&  0.0 \\
500&  44.5&&  25.8&&  25.4&&  13.5&  0.0&  0.0&  0.0 \\
	\bottomrule 
	\multicolumn{1}{r}{\vspace{-2ex}}\\
	\multicolumn{10}{l}{\footnotesize$^\dagger$Average excludes one instance that exceeded an hour in computation time.} \\
	\multicolumn{10}{l}{\footnotesize$^{\dagger\dagger}$Average excludes two instances that exceeded an hour.} \\
	\multicolumn{10}{l}{\footnotesize$^*$All three instances exceeded an hour.} 
\end{tabular}
\end{table}

As one might expect, the integer L-shaped implementations are faster than solving the deterministic equivalent MILP, because they exploit the scenario-based block structure of two-stage stochastic programs.  We also see that the integer L-shaped method can be significantly accelerated by solving the exact subproblem with CP rather than MILP (to obtain upper bounds and generate the integer cut), since CP is more effective for this type of scheduling problem. 

It is clear from Table~\ref{table:CompPerf1} that all four implementations of LBBD substantially outperform the integer L-shaped method, even when the latter uses CP.  Furthermore, the two branch-and-check implementations scale much better than standard LBBD, due mainly to time spent in solving the master problem in standard LBBD. This confirms the rule of thumb that branch and check is superior when solving the master problem takes significantly longer than solving the subproblems.  The results also indicate that analytical Benders cuts are more effective than unstrengthened nogood cuts in both standard LBBD and branch and check.  

In summary, branch and check with analytical Benders cuts is the best of the seven methods for these test instances.  In particular, it is far superior to the integer L-shaped method, as it easily solves instances with 500 scenarios, while the integer L-shaped method cannot deal with more than 10 scenarios within an hour of computation time.

Table~\ref{table:AlgAnalysis} provides a more detailed comparison of the integer \mbox{L-shaped} method with the branch-and-check implementations.  The L-shaped method with CP is shown, as we have seen that it is faster than solving the subproblem with MILP.  Interestingly, solving a CP  formulation of the subproblem is much faster than solving the LP relaxation of an MILP formulation.  This illustrates the computational cost of using the larger MILP formulation.  We also see that the stronger analytical cuts substantially reduce the number of times the subproblem must be solved, and therefore the number of cuts generated and the resulting size of the master problem.  Furthermore, the number of subproblem calls is roughly constant as the number of scenarios increases.  Finally, the subproblem solutions consume about half of the total computation time in the branch-and-cut algorithms.  Previous experience suggests that for best results, the computation time should, in fact, be about equally split between the master problem and subproblem (\citeauthor{CirCobHoo16} \citeyear{CirCobHoo16}).

\begin{table}[!t]
\small
\setlength{\tabcolsep}{4pt} 
\centering
\caption{Analysis of the integer L-shaped method with CP subproblems and two branch-and-check algorithms.  Each number is an average over 3 problem instances.}
\label{table:AlgAnalysis}
\vspace{3ex}
\begin{tabular}{r | r@{}r@{}@{\hspace{2.2ex}}r@{}r@{}@{\hspace{1.5ex}}r@{}r@{}@{\hspace{2ex}}r@{}r@{}@{\hspace{0ex}}r@{}r | rrr@{\hspace{1ex}}r | rr@{\hspace{4ex}}r@{\hspace{0.5ex}}r }
	\multicolumn{1}{c}{} & \multicolumn{10}{c}{Integer L-shaped with CP}
	& \multicolumn{4}{c}{B\&C with nogood cuts} 
	& \multicolumn{4}{c}{B\&C with analytical cuts}\\
	\cmidrule{2-19} 
	\multicolumn{1}{c}{} & \multicolumn{6}{c}{Time (sec)} & \multicolumn{4}{c}{Statistics}
	& \multicolumn{2}{|c}{Time (sec)} & \multicolumn{2}{c}{Statistics}
	& \multicolumn{2}{|c}{Time (sec)} & \multicolumn{2}{c}{Statistics}\\
	\cmidrule{2-19} 
	Scenarios & 
	Total && \multicolumn{1}{@{\hspace{-3ex}}c@{\hspace{-1ex}}}{CPsub} && \multicolumn{1}{@{\hspace{-2ex}}c@{\hspace{-1ex}}}{LPsub} && Cuts && Calls
	&& Total  & \multicolumn{1}{@{\hspace{-1ex}}c@{\hspace{-1ex}}}{CPsub} & Cuts & Calls
	& Total & \multicolumn{1}{@{\hspace{0ex}}c@{\hspace{-1ex}}}{CPsub} & Cuts & Calls \\
	\toprule
1&  27.9&&  0.9&&  2.1&&  450&&  452&&  1.8&  0.5&  282&  150&  1.6&  0.2&  82&  47 \\
5&  149.3&&  5.6&&  16.4&&  2692&&  541&&  3.3&  1.7&  1289&  144&  1.7&  0.5&  311&  40 \\
10&  437.8&&  15.8&&  73.4&&  5114&&  515&&  5.1&  2.6&  2390&  134&  2.8&  1.0&  702&  46 \\
50&  2517.8&$^{\dagger}$ &  97.3&$^{\dagger}$\hspace{-2.75ex}&  500.2&$^{\dagger}$\hspace{-2.75ex}&  20002&$^{\dagger}$\hspace{-2.75ex}&  401&$^{\dagger}$\hspace{-2.75ex}&  33.3&  16.9&  12616&  148&  18.0&  7.6&  3906&  51 \\
100&  *&&  *&&  *&&  *&&  *&&  80.1&  42.0&  25880&  152&  30.2&  11.3&  7607&  50 \\
500&  *&&  *&&  *&&  *&&  *&&  416.2&  187.6&  127404&  150&  166.2&  55.7&  35409&  47 \\
	\bottomrule
	\multicolumn{1}{r}{\vspace{-2ex}}\\
	\multicolumn{19}{l}{\footnotesize$^{\dagger}$Average excludes two instances that exceeded an hour.} \\
	\multicolumn{19}{l}{\footnotesize$^*$Computation terminated for all 3 instances after one hour.} 
\end{tabular}
\end{table}

Given the computational burden of solving the LP relaxation of the MILP subproblem, we experimented with running the integer \mbox{L-shaped} method with only integer cuts.  This obviates the necessity of solving the LP relaxation of an MILP model.  We also solved instances with 14 and 18 as well as 10 tasks and with 4 facilities as well as 2.  The results appear in Table~\ref{table:CompPerf2}. The three implementations shown in the table are exactly the same except for the cuts used and therefore permit a direct comparison of the effectiveness of the cuts.  Comparison with Tables~\ref{table:CompPerf1} and~\ref{table:AlgAnalysis} reveals that the integer \mbox{L-shaped} method actually runs faster using only integer cuts, without any classical Benders cuts obtained from the LP relaxation.  We also see that the analytical cuts are more effective than unstrengthened nogood cuts in nearly every instance, and much more effective than integer cuts, which are quite weak.  Finally, branch and cut is increasingly superior to even this accelerated version of the integer L-shaped method as the instances scale up, and therefore far superior to the standard method.

\begin{table}[!t]
\small
\centering
\caption{Average computation time in seconds over 3 instances for various types of cuts.}
\label{table:CompPerf2}
\vspace{3ex}
\begin{tabular}{rr | rrr | rrr }
	\multicolumn{2}{c}{} & \multicolumn{3}{c}{2 facilities} & \multicolumn{3}{c}{4 facilities} \\
	\cmidrule{3-5} \cmidrule{6-8}
	&		&L-shaped
	                &B\&C	&B\&C	&L-shaped
	                                        &B\&C	&B\&C  \\
Tasks
    &Scenarios
            &integer&nogood	&analytic
                                    &integer&nogood	&analytic \\
	&		&cuts only&cuts	&cuts	&cuts only	
	                                        &cuts	&cuts \\												
	\toprule 
	10		&1		&1.8	&1.2	&0.3	&78.1	&0.5	&0.4		\\
	&5		&8.7	&3.0	&1.2	&906.5	&2.5	&2.3		\\
	&10		&15.3	&6.0	&2.3	&2213.0$^\dagger$ \hspace{-1.75ex}
	&3.9	&4.6		\\
	&50		&90.9	&27.4	&10.7	&*		&24.4	&19.6		\\
	&100	&217.9	&50.0	&22.0	&*		&41.9	&33.5		\\
	&500	&1318.3	&329.3	&111.0	&*		&268.3	&205.2		\\
	\midrule
	14	&1		&48.9	&4.3	&1.5	&2403.9$^{\dagger\dagger}$ \hspace{-2.75ex}
	&0.9	&0.9		\\
	&5		&229.5	&16.3	&5.7	&2402.8$^{\dagger\dagger}$ \hspace{-2.75ex}
	&5.6	&3.1		\\
	&10		&284.7	&37.7	&9.0	&*		&12.4	&5.9		\\
	&50		&1850.6	&186.0	&31.6	&*		&88.0	&41.9		\\
	&100	&2810.4$^{\dagger\dagger}$ \hspace{-2.75ex} 
	&411.0	&70.2	&*		&189.5	&75.2		\\
	&500	&*		&2431.9$^\dagger$ \hspace{-1.75ex} 
	&494.8	&*		&854.7	&504.9		\\
	\midrule
	18	&1		&1358.6$^\dagger$ \hspace{-1.75ex} 
	&208.3	&15.2	&1346.8$^\dagger$ \hspace{-1.75ex}
	&1.5	&1.3		\\
	&5		&229.5	&16.3	&5.7	&2402.8$^{\dagger\dagger}$ \hspace{-2.75ex}
	&5.6	&3.1		\\
	&10		&3477.2$^\dagger$ \hspace{-1.75ex}	
	&2184.0	&113.4	&*		&116.4	&35.2		\\
	&50		&* 		&*		&1138.1	&*		&458.2	&148.5		\\
	&100	&*		&*		&2298.7	&*		&943.7	&285.6		\\
	&500	&*		&*		&*		&*		&2804.5$^{\dagger\dagger}$	\hspace{-2.75ex}
	&1318.9		\\
	\bottomrule 
	\multicolumn{1}{r}{\vspace{-2ex}}\\
	\multicolumn{8}{l}{\footnotesize$^\dagger$Average excludes one instance that exceeded an hour in computation time.} \\
	\multicolumn{8}{l}{\footnotesize$^{\dagger\dagger}$Average excludes two instances that exceeded an hour.} \\
	\multicolumn{8}{l}{\footnotesize$^*$All three instances exceeded an hour.} 
\end{tabular}
\end{table}

\section{Conclusion}
\label{Sec:Conclusion}

In this study, we applied logic-based Benders decomposition (LBBD) to two-stage stochastic optimization with a scheduling task in the second stage.  While Benders decomposition is often applied to such problems, notably in the integer L-shaped method, the necessity of generating classical Benders cuts requires that the subproblem be formulated as a mixed integer/linear programming problem and cuts generated from its continuous relaxation.  We observed that this process incurs substantial computational overhead that LBBD avoids by generating logic-based cuts directly from a constraint programming model of the scheduling subproblem.  Although the integer cuts used with the L-shaped method can be regarded as a special case of logic-based Benders cuts, they are extremely weak, even weaker than simple nogood cuts often used in an LBBD context.  Furthermore, the type of subproblem analysis that has been used for past applications of LBBD permits much stronger logic-based cuts to be derived, again without the overhead of obtaining a continuous relaxation.  

Computational experiments found that, due to these factors, LBBD solves a generic stochastic planning and scheduling problem much more rapidly than the integer L-shaped method.  The speedup is several orders of magnitude when a branch-and-check variant of LBBD is used.  This outcome suggests that LBBD could be a promising approach to other two-stage stochastic and robust optimization problems with integer or combinatorial recourse, particularly when the subproblem is relatively difficult to model as an integer programming problem.

\clearpage






\end{document}